\documentclass[12pt,notitlepage]{amsart} 
\usepackage{amsthm,amsmath,amssymb,array,delarray,epsfig,latexsym} 

\newlength{\sh}
\newlength{\jim} 
\newlength{\smale} 
\newlength{\jmr}
\newlength{\jfc}
\newlength{\bernd}
\newlength{\ioan}
\newlength{\agk}
\newtheorem{lemma}{Lemma}
\newtheorem{prop}{Proposition}
\newtheorem{disc}{Univariate Discriminant Complexity Problem} 
 
\newtheorem{proba}{Problem A} 
 
\newtheorem{probb}{Problem B} 
 
\newtheorem{dfn}{Definition}

\newtheorem{main}{Main Theorem}

\newtheorem{realsmale}{Real Analogue of Smale's 17$^\thth$ Problem}

\newtheorem{thm}[main]{Theorem}
\newtheorem{rem}{Remark}	
\newtheorem{ex}{Example}

\newcommand{\eps}{\varepsilon}

\newcommand{\cN}{\mathcal{N}}
\newcommand{\cO}{O} 
\newcommand{\cF}{\mathcal{F}}

\newcommand{\np}{\mathbf{NP}}
\newcommand{\pp}{\mathbf{P}}
\newcommand{\bpp}{\mathbf{BPP}}

\newcommand{\thth}{{\underline{\mathrm{th}}}}

\newcommand{\Q}{\mathbb{Q}}

\newcommand{\R}{\mathbb{R}}
\newcommand{\C}{\mathbb{C}}
\newcommand{\N}{\mathbb{N}}
\newcommand{\Z}{\mathbb{Z}}

\newcommand{\Zn}{\Z^n}

\newcommand{\Rn}{\R^n}
\newcommand{\balpha}{{\bar{\alpha}^*}}

\newcommand{\Cn}{\C^n}

\renewcommand{\qed}{$\blacksquare$}
\newcommand{\dia}{$\diamond$}

\newcommand{\cC}{\mathcal{C}}
\newcommand{\cD}{\mathcal{D}}
\newcommand{\cK}{\mathcal{K}}

\begin{document}  

\title[On Solving Fewnomials in Fewnomials Time]{On Solving 
Fewnomials Over Intervals in Fewnomial Time} 

\author{J. Maurice Rojas$^\star$}\thanks{
\mbox{}\hspace{-.7cm}
$\star$ Partially funded by Hong Kong/France PROCORE Grant 
\#9050140-730, a Texas A\&M University Faculty of Science 
Grant, and NSF Grants DMS-0138446 and DMS-0211458. }

\author{Yinyu Ye$^\dagger$}\thanks{
\mbox{}\hspace{-.55cm}
$\dagger$ Partially funded by NSF Grant 
DMS-9703490 and the City University of Hong Kong. The second author's work 
was done while visiting the Department of Mathematics at the City University 
of Hong Kong.} 

\address{Department of Mathematics, Texas A\&M University, TAMU 3368, 
College Station, TX \ 77843, USA } 

\email{rojas@math.tamu.edu\\ {\it Web-Page:} 
http://www.math.tamu.edu/\~{}rojas}

\address{Department of Management Science and Engineering, School of 
Engineering, Stanford University, Stanford, CA \ 94305, USA
} 

\email{yinyu-ye@stanford.edu\\ {\it Web-Page:} 
http://www.stanford.edu/\~{}yyye} 

\date{\today} 

\begin{abstract} 
Let $f$ be a degree $D$ univariate polynomial with real coefficients 
and exactly $m$ monomial terms. We show that in the special case 
$m\!=\!3$ we can approximate within $\eps$ all the roots of $f$ in 
the interval $[0,R]$ using just $O\!\left(\log(D)\log\left(D\log 
\frac{R}{\eps}\right)\right)$ arithmetic operations. In particular, 
we can count the number of roots in any bounded interval using just 
$O(\log^2 D)$ 
arithmetic operations. Our speed-ups are significant and near-optimal: The 
asymptotically sharpest previous complexity upper bounds for 
both problems were super-linear in $D$, while our algorithm 
has complexity close to the respective complexity lower bounds. We also 
discuss conditions 
under which our algorithms can be extended to general $m$, and a connection 
to a real analogue of Smale's 17$^\thth$ Problem.  
\end{abstract} 

\mbox{}\\
\vspace{-.9in} 
\maketitle 

\section{Introduction}
\label{sec:intro}
Real-solving --- the study of solving systems of polynomial equations over 
the real numbers --- occupies a curious position within computational algebraic 
geometry. From the point of view of computational complexity, classical 
algebraic geometry has left real-solving almost completely untouched. 
For example, rigorous lower bounds for the arithmetic complexity of finding 
approximations to the roots of polynomial systems didn't appear until the 
work of Renegar in the late 1980's \cite{renuni,renegar}, 
and finding optimal bounds continues to be the subject of much active 
research \cite{bshouty,bcss,mp99,mygcp,gls01,mpr03}. 

As for counting real roots, there is a beautiful result that one can bound  
their number independently of the degrees of the underlying polynomials: 
Askold Khovanski proved an explicit upper bound singly exponential in the 
number of variables and the total number of monomial terms \cite{sparse}. 
While this {\bf fewnomial} bound is far from optimal in higher dimensions 
\cite{tri}, it is significantly smaller than the number of complex roots 
when the underlying polynomials have sufficiently high degree. 
(Sparse polynomials are sometimes also known as lacunary polynomials and, 
over $\R$, are a special case of fewnomials --- a more general class 
of analytic functions of parametrized complexity \cite{sparse}.)  

We are then naturally lead to suspect that a similar improvement is possible 
for the harder problem of approximating the real roots. So can one solve 
{\bf sparse} polynomial systems 
over the real numbers significantly faster than via the usual algorithms 
based on complex algebraic geometry? The existence of general speed-ups 
of this nature is still an open problem, even in the univariate case: 
For example, until the present paper, it was still unknown whether the real 
roots of a univariate trinomial of degree $D$ could be approximated within a 
number of arithmetic operations sub-linear in $D$ \cite{mp99,mygcp,
gls01,mpr03}. 

We now answer this last question affirmatively as follows: 
Let $I\!\subseteq\!\C$ be any subset, $f\!\in\!\C[x_1]$ a degree $D$ 
polynomial, and suppose $\zeta_1,\ldots,\zeta_m$ are all the distinct 
roots of $f$ in the region $I$. By {\bf $\pmb{\eps}$-approximating the roots 
of $\pmb{f}$ in $\pmb{I}$} we will mean finding complex numbers 
$z_1,\ldots,z_m\!\in\!I$ such that for any $i\!\in\!\{1,\ldots,m\}$, we have  
$|\zeta_i-z_j|\!<\!\eps$ for some $j\!\in\!\{1,\ldots,m\}$. 
(So $\eps$-approximating in $I$ implies a correct root count in  
$I$ as well.) When $I\!\subseteq\!\R$ and $f\!\in\!\R[x_1]$, we will further 
stipulate that the $z_i$ all lie in $\R$. 
\begin{thm} 
\label{thm:uni}
Let $R,\eps\!>\!0$ and suppose $f\!\in\!\R[x_1]\setminus\{0\}$ has degree 
$D$ and at most $3$ monomial terms. Then we can $\eps$-approximate all the 
roots of $f$ in the closed interval $[0,R]$ using just 
$\cO\!\left(\log(D)\log\left(D\log \frac{R}{\eps}\right)\right)$ 
arithmetic operations. In particular, we can count exactly the number of roots 
of $f$ in any bounded interval using just $\cO(\log^2 D)$ arithmetic 
operations. 
\end{thm}

\noindent
As is standard, we count arithmetic operations as field 
operations in the field over $\Q$ generated by the coefficients of $f$. 
Note also that our underlying algorithm handles degenerate roots 
(i.e., roots of multiplicity $>\!1$) with no difficulty. 
\begin{rem}
Throughout this paper, all $O$-constants and $\Omega$-constants are
absolute and effectively computable. \dia
\end{rem}
\begin{rem} 
Our speed-ups are significant: The current asymptotically sharpest 
(sequential worst-case) arithmetic 
\mbox{complexity} upper bound for $\eps$-approximating 
the roots of a general degree $D$ univariate polynomial in the open disc 
$\{z \; | \; |z|\!\leq\!R\}$ is 
$\cO\!\left(D\log^5(D)\log\log\frac{R}{\eps}\right)$ \cite{binipan,neffreif}.  
The analogous upper bound for counting real roots in an 
interval is \scalebox{.8}[1]{$\cO\!\left(D\log^2(D)\log\log D\right)$}, 
via the technique of  {\bf Sylvester-Habicht sequences} \cite{marie,lickroy}. 
In particular, no sharper upper bounds were known before for the 
problems considered in Theorem \ref{thm:uni} above.  \dia 
\end{rem} 
\begin{rem}  
Our speed-ups are also near-optimal: For any fixed $D,R\!\geq\!2$, our 
$\eps$-approximation algorithm matches (up to an asymptotically constant 
multiple) the $\Omega\!\left(\log\log \frac{1}{\eps}\right)$ arithmetic 
complexity lower bound known for $\eps$-approximating the roots of $x^2-N$ 
where $1\!\leq\!N\!\leq\!2$ \cite{bshouty}.  
As for counting the roots, there do not appear to be any explicit 
lower bounds known for the arithmetic complexity of counting the real 
roots of $m$-nomials. However, one should note that the arithmetic complexity 
of just evaluating a degree $D$ monomial is $\Omega(\log D)$ in the worst 
case \cite{fux,gugu}. \dia  
\end{rem} 
Our algorithms are based on an earlier hybrid algorithm of the
second author which combines bisection and Newton iteration, a
new observation on the {\bf Sturm sequences} of 
trinomials (Theorem \ref{thm:bisturm} of Section \ref{sub:count}), 
and some analytic estimates (Theorem \ref{thm:global} of Section 
\ref{sec:back}). We in fact give a more general algorithm,  
in Theorem \ref{thm:big} of Section \ref{sec:back},  
that applies to certain univariate $m$-nomials. However, for $m\!\geq\!4$,  
there are two main obstructions to showing that our more general algorithm 
has complexity sub-linear in $D$. We refer to these obstructions as 
Problems A and B, and describe the first now. 
\begin{proba} Is there an absolute contant $\kappa$ such that one can count 
exactly the number of real roots (in any input interval) of arbitrary 
$m$-nomials of degree $D$ using, say, just $O\!\left(\log^{\kappa m} D\right)$ 
arithmetic operations?  \dia 
\end{proba} 
One may note that we have hedged our bets in Problem A by 
letting the complexity bound increase exponentially in $m$. 
This is motivated by the $\mathbf{NP}$-hardness of the multivariate analogue 
of Problem A.
\begin{prop} 
Suppose $f\!\in\!\Z[x_1,\ldots,x_n]$ is an $m$-nomial and 
we measure the {\bf size} of $f$ as the total number of 
bits necessary to write the binary expansions of all the coefficients 
and exponents of $f$. Then it is $\mathbf{NP}$-hard (in the classical 
Turing sense, using the preceding notion of size) to decide whether $f$ has a 
real root or not. \qed 
\end{prop} 

\noindent
(The proof follows immediately from a standard reduction of the special case 
of $n$-variate quartic polynomials to $\mathbf{SAT}$ \cite{bcs}.) 
We present some examples in Section \ref{sub:count} revealing that 
the classical approach of Sturm sequences will most likely need to be 
abandonded for $m\!>\!3$. Nevertheless, speed-ups for $m$ {\bf fixed} still 
appear possible and are of considerable practical interest for large $D$.  

Noting that $\eps$-approximation extends naturally to polynomial 
systems (by $\eps$-approximating each coordinate separately), a consequence of 
our univariate trinomial algorithm is the following result which may be of use 
for solving general pairs of bivariate trinomials. 
\begin{thm} 
\label{thm:tri} 
Suppose $f\!\in\!\R[x_1,x_2]$ has exactly $3$ monomial terms, 
$D$ is the degree of $f$, and $\delta\!>\!0$ is any constant. Then, using 
just 
\mbox{$\cO\!\left(\log(D)\log\left(D\log\frac{R}{\eps}\right)\right)$} 
arithmetic operations and $\cO\left(\log^{2+\delta}D\right)$ bit operations, 
we can $\eps$-approximate all isolated inflection points, vertical 
tangents, and singular points of the curve $\left\{(x_1,x_2)\!\in\!
[0,R]^2 \; | \; f(x_1,x_2)\!=\!0\right\}$.  
\end{thm}
\noindent 
The asymptotically sharpest (sequential worst-case) 
arithmetic complexity upper bound for $\eps$-approximating the 
roots of a system of $n$ polynomial equations in $n$ variables are at 
best polynomial in 
$\log\log \frac{R}{\eps}$ and a quantity sometimes attaining $D^n$ 
(see, e.g., \cite{mp99,mygcp,gls01,mpr03}). The importance of Theorem 
\ref{thm:tri} lies in that it is a first step toward a practical 
two-dimensional analogue of bisection. The latter in turn is a first step 
toward generalizing our algorithms here to higher dimensions. 

We leave the bit complexity of our algorithms for a future paper. In 
particular, unless stated otherwise, our underlying model of computation 
here is the BSS model (with 
inequality) over $\R$ \cite{bcss}. Those unfamiliar with this model can simply 
think of such a machine as an ordinary Turing machine augmented with registers 
that allow arithmetic and inequalities with real numbers as well as bits. 
However, let us touch upon a deeper question: The connection of our work to 
Smale's notion of {\bf approximate roots} (see, e.g., \cite{smale,bcss} and 
Section \ref{sec:back} below). 
\begin{dfn} 
Suppose $f_1,\ldots,f_n\!\in\!\R[x_1,\ldots,x_n]$, $F\!:=\!(f_1,\ldots,f_n)$, 
and the total number of distinct exponent vectors in the monomial 
term expansions of $f_1,\ldots,f_n$ is $m$. We then call 
$F$ a {\bf real $\pmb{m}$-sparse $\pmb{n\times n}$ polynomial system}. \dia 
\end{dfn} 
\begin{dfn} 
Let $B\!\subseteq\!\Cn$ be any open ball and 
$F : B\longrightarrow \Cn$ an analytic function. Also, 
given any $z_0\!\in\!B$, let us define the sequence 
$(z_i)^\infty_{i=0}$ by $z_i\!:=\!z_{i-1}-\mathrm{Jac}(F)|^{-1}_{z_{i-1}}
F(z_{i-1})$ for all $i\!\in\!\N$. Then, if there is a root $\zeta\!\in\!B$ 
of $F$ 
such that $|z_{i+1}-\zeta|\!\leq\!8(\frac{1}{2})^{2^i}|z_0-\zeta|$ for all $i$
(i.e., if Newton iteration for $F$, starting at $z_0$, converges quadratically
to a root of $F$), we call $z_0$ an {\bf approximate root} of $F$ with {\bf
associated root} $\zeta$. \dia
\end{dfn}
\begin{realsmale} 
For fixed $n$, can all real roots of an $m$-sparse $n\times n$ 
polynomial system be found approximately, on the average, 
in polynomial time with a uniform algorithm? More precisely, 
let $F$ be an $m$-sparse $n\times n$ polynomial system with 
maximal exponent $D$ and coefficients that are, say, independent 
standard real Gaussian random variables. Is 
there a uniform algorithm that finds a set of 
approximate roots close to all the real roots of $F$, with 
average-case arithmetic complexity $O((m\log D)^{\nu})$ for some constant 
$\nu$ depending only on $n$? \dia 
\end{realsmale} 

\noindent
The polynomial system $\left((x_1-1)(x_1-2),\ldots,(x_n-1)(x_n-2)\right)$ 
clearly shows that fixing $n$ is necessary in our real analogue above. 
The appellation {\bf uniform} merely 
emphasizes that there be a single algorithm, with all steps explicit and
constructive, which works for all inputs \cite{bcss}.

The original statement of Smale's 17$^\thth$ Problem \cite{next} (see also 
\cite[Pg.\ 287]{21})  
differs from our analogue above as follows: (1) $n$ is allowed to 
vary (so one seeks an {\bf absolute} constant $\nu$), (2) one 
averages over choices of complex coefficients, 
and (3) one instead asks for a {\bf single} complex approximate root. 
Smale also left the underlying probability distribution unspecified. 
We observe that a positive answer to the following variant of Problem A would 
be quite useful in the direction of our real analogue: 

\noindent 
{\bf Stochastic Version of Problem A:} {\em Suppose 
$f$ is a univariate $m$-nomial of degree $\leq\!D$ with coefficients that are 
real standard Gaussian random variables. Is there a uniform algorithm that 
counts exactly the number of roots of $f$ (in any input interval), with 
{\bf average-case} arithmetic complexity 
$O((m\log D)^{\kappa})$ for some absolute constant $\kappa$? \dia  } 

Smale's original 17$^\thth$ Problem remains unsolved, 
although a partial affirmative answer (an algorithm containing 
a non-constructive step which may be called many times) was found by Shub 
and Smale in the mid-1990's \cite{ss5,bcss}. 

\subsection{Earlier Work on Approximating Roots}\mbox{}\\
\label{sub:early} 
When can we solve a polynomial system in time polylogarithmic in 
the degree of the underlying complex algebraic set? 
An affirmative answer is trivial for the special case of a single polynomial 
with $\leq\!1$ monomial term, and was known at least since the 
mid-1970's for the (univariate) binomial case, e.g., \cite{brent} and 
\cite[Sec.\ 4]{ye}. On the other hand, little 
seems to be known about the case of $3$ or more monomial terms: To the best of 
the authors' knowledge, the only result close in spirit to Theorem 
\ref{thm:uni} is a result of Daniel Richardson \cite{richardson} implying that 
the arithmetic complexity of counting the number of real roots of 
$c_1+c_2x^d+c_3x^D$ (with $0\!<\!d\!<\!D$) is polynomial in $d$ and $\log D$.   

Not much else seems to be known about the intrinsic complexity of real solving, 
or even real root counting, for univariate polynomials with  
$4$ or more monomial terms. One known speed-up over the usual univariate 
algorithms over $\C$ is \cite[Main Theorem 1.2]{real} which gives an 
arithmetic complexity bound of $\cO\!\left(m\log(D)\log\left(D\log
\frac{R}{\eps}\right)\right)$ for $\eps$-approximating the roots of 
\[c_1x^{a_1}+\cdots+c_kx^{a_k}-c_{k+1} 
x^{a_{k+1}}-\cdots-c_mx^{a_m},\] when $0\!\leq\!a_1\leq \cdots 
\leq\!a_m\!=\!D$ are integers and all the $c_i$ are positive real numbers. 

It is also not difficult via the results above to construct 
various systems of multivariate trinomials which admit super-fast 
solving in the sense of Theorem \ref{thm:uni} as well. For example, 
if one has $f_1\!\in\!\R[x_1]$, $f_2\!\in\!\R[x_1,x_2],\ldots$ 
$f_n\!\in\!\R[x_1,\ldots,x_n]$, and all the $f_i$ are trinomials of 
degree $\leq\!D$, then 
to solve $F\!:=\!(f_1,\ldots,f_n)$ one can simply solve $f_1$ first and then 
recursively solve the resulting smaller system. Letting $\cN_F$ 
denote the number of roots in the nonnegative orthant of such an $F$, 
it is easily checked that $\cN_F\!\leq\!2^n$ (see, e.g., \cite[Thm.\ 3, 
Part (c)]{tri}). We can then easily derive, 
via Theorem \ref{thm:uni},  an arithmetic complexity upper bound of 
$\cO\!\left(\cN_F\log^2(D)\log\log\frac{R}{\eps}\right)$ 
for the more general problem where we instead $\eps$-approximate all the 
roots of $F$ in the {\bf orthant-wedge} 
\[W^n_R\!:=\!\left\{(x_1,\ldots,x_n)\!\in\!\Rn \; \left| 
\; x^2_1+\cdots+x^2_n\!\leq\!R^2 \ 
, \ x_i\!\geq\!0 \text{ \ for \ all \ } i\right. \right\}.\] 
Extending these results to general trinomial systems, not to 
mention general sparse systems, remains an open problem. 
Nevertheless, for a general {\bf bi}nomial system $F$ with exactly 
$\cD\!<\!\infty$ roots in $\Cn$ (counting multiplicities), 
one can $\eps$-approximate all its roots in $W^n_\R$ using 
just \scalebox{.8}[1]{$O\!\left((\log \cD)\left(B^3\log^2(n)+\log\log 
\frac{R}{\eps}\right) \right)$} arithmetic operations and 
\scalebox{1}[1]{$O(n^3B^{2+\delta} \log^{2+\delta}n)$}  
bit operations, where $B$ is the total number of bits needed to write 
down the exponents of $F$ and $\delta\!>\!0$ is an arbitrary constant  
\cite[Main Theorem 1.3]{real}. (The bound in \cite{real} is 
written in a slightly different manner and we have here taken the liberty of 
improving the bit complexity portion by employing a recent result 
of van der Kallen on computing the Hermite normal form of an 
integral matrix \cite{vanderk}.)  

\subsection{A New Algebraic Observation and Earlier Work on Counting Roots} 
\mbox{}\\
\label{sub:count}
As for merely counting the real roots, there are still surprising gaps in our 
knowledge. For instance, the main general algebraic techniques for real 
root counting --- Sturm-Habicht sequences \cite{sturm,habicht,marie,lickroy} 
and Hermite's Method \cite{hermite,prs,marie} ---  originated in the 
19$^\thth$ century and are still hard to improve upon. 
\begin{dfn}
\label{dfn:sturm}
For any sequence of real numbers $s\!:=\!(s_1,\ldots,s_k)$, the {\bf number of
sign alternations of $\pmb{s}$}, $\pmb{N_s}$, is the number of pairs
$(j,j')$ with $1\!\leq\!j<j'\!\leq\!k$, $s_js_{j'}\!<\!0$, and $s_i\!=\!0$
for all $i\!\in\!\{j+1,\ldots,j'-1\}$. Also, for any polynomial
$f\!\in\!\R[x_1]$,
define the sequence $(p_i)^\infty_{i=0}$ where
$p_0\!:=\!f$, $p_1\!:=\!f'$, $p_{i+2}\!:=\!q_ip_{i+1}-p_i$, $q_i$
is the quotient of $\frac{p_i}{p_{i+1}}$, and $p_K$ is the last element of
the sequence not identically equal to $0$. We then call $(p_0,\ldots,p_K)$ 
the {\bf Sturm sequence} of $f$. \dia
\end{dfn}
\begin{lemma}
\cite{sturm,marie}
\label{lemma:sturm}
Following the notation of Definition \ref{dfn:sturm}, let 
$a,b\!\in\!\R$
with $a\!\leq\!b$, $A\!:=\!(p_0(a),\ldots,p_K(a))$,
and $B\!:=\!(p_0(b),\ldots,p_K(b))$. Then the number of distinct roots 
of $f$ in the open interval $(a,b)$ is exactly $N_A-N_B$. In 
particular,
[$a_1<\cdots<a_m$ and $c_1,\ldots,c_m\!\in\!\R$]
$\Longrightarrow c_1x^{a_1}+\cdots
+c_mx^{a_m}$ has no more than $N_{(c_1,\ldots,c_m)}$ ($\leq\!m-1$)
positive \mbox{roots. \qed}
\end{lemma}
\begin{rem}
The latter part of the above lemma is {\bf
Descartes' Rule of Signs} which dates back to
the famous philosopher's 1637 book La Geometrie (see also
\cite[Pg.\ 160]{descartes}). \dia
\end{rem}
                                                                                
The following observation on the Sturm sequence of a trinomial,
which we use in proving Theorem \ref{thm:uni}, may be of independent
interest.
\begin{thm}
\label{thm:bisturm}
Following the notation above, suppose $f$ has at most $3$ monomial terms.
Then $K\!\leq\!3\lceil\log_2 D\rceil +2$ and $(p_1,\ldots,p_K)$ consists
solely of binomials, monomials, and/or constants. In particular, the entire
Sturm sequence of $f$ can be evaluated at any real number using just
$\cO(\log^2 D)$ arithmetic operations.
\end{thm}

Extending our last theorem to polynomials with more monomials appears 
unlikely: First, note that the quotient of the
division of two binomials can be quite non-sparse, e.g.,
$\frac{x^{N+1}-1}{x-1}\!=\!x^N+\cdots+x+1$. So the expansion of any
intermediate quotients must be avoided. Furthermore, the {\bf tetra}nomial case
already gives some indication that Sturm sequences may in fact have to be
completely abandoned: the following example shows that, for all
$D\!>\!2$, the fourth element of the Sturm sequence of a degree $2D$
tetranomial can already be a $(D+1)$-nomial.
\begin{ex}
Consider the tetranomial $p_0(x)\!:=\!x^{2D}+x^{D+1}+x^D+1$ and its
derivative
\mbox{$p_1\!:=\!p'_0\!=\!2Dx^{2D-1}+(D+1)x^D+Dx^{D-1}$.} The resulting Sturm
sequence then continues with $p_2\!=\!-\frac{D-1}{2D}
x^{D+1}-\frac{1}{2}x^D-1$ and from here it is easy to see (by a
writing a simple recursion for the resulting long division) that
the quotient $q_3$ of $\frac{p_1}{p_2}$ is a polynomial of degree $D-2$
with exactly $D-1$ monomial terms. Thus, $p_3\!:=\!q_3p_2-p_1$ has
degree $D$ and at least $D+1$ monomial terms. \dia
\end{ex}
\noindent
We are also willing to conjecture that the maximal length of the Sturm
sequence of a degree $2D$ tetranomial is $\Omega(D)$. ({\tt Maple} experiments
have verified this up to $D\!=\!150$.) Nevertheless, while the behavior of
tetranomials is thus more complicated, this example need not rule out a more
clever method to circumvent these difficulties.

\subsection{A Connection to Discriminants}\mbox{}\\  
If one insists on relying on Sturm sequences then one is
naturally lead to the univariate sparse discriminant. Briefly,
given an $n$-variate $m$-nomial $f$ with indeterminate coefficient vector $C$ 
and set of exponents $A$, its {\bf sparse discriminant} (or 
{\bf $\pmb{A}$-discriminant}), $\Delta_A(f)$, is the
unique (up to sign) irreducible polynomial in $\Z[C]\setminus\{0\}$ of
lowest degree which vanishes whenever $C$ is specialized so that $f$
has a root in $(\C\!\setminus\!\{0\})^n$ in common with 
$\frac{\partial f}{\partial x_1},\ldots, \frac{\partial f}{\partial x_n}$ 
\cite{gkz94}. It is not hard to see that if 
(a) $n\!=\!1$, (b) the coefficients of $f$ are all constants, 
and (c) $A$ is exactly the set of exponent vectors of $f$, 
then we have $\Delta_A(f)\!=\!0 \Longleftrightarrow 
K\!=\!1$ in the Sturm sequence $(p_0,\ldots,p_K)$ of $f$. So the following is 
clear. 
\begin{prop}
\label{prop:disc}
Suppose $f\!\in\!\R[x_1]$ has set of exponent vectors $A$ 
and one can compute the signs of the Sturm sequence of $f$ evaluated 
at any $r\!\in\!\R$ using just $T(A)$ arithmetic operations. Then one can 
decide the vanishing of $\Delta_A(f)$ using just $T(A)$ arithmetic 
operations. \qed
\end{prop}

\noindent 
It would be quite enlightening to understand the converse of Proposition 
\ref{prop:disc} and its possible obstructions. 
                                                                                
Perhaps more than coincidentally  --- recalling Theorem \ref{thm:uni} --- 
$\Delta_A(f)$ {\bf can} be computed in
polynomial time when $f$ is a trinomial: up to sign, the formula is
simply \[\Delta_{\{0,a_2,a_3\}}(c_1+c_2x^{a_2}+c_3x^{a_3})\!=\!
a^{a_3}_3c^{a_2}_3c^{a_3-a_2}_1+(-1)^{a_3-1}a^{a_2}_2(a_3-a_2)^{a_3-a_2}
c^{a_3}_2\] when $a_2$ and $a_3$ are relative prime and $0\!<\!a_2\!<\!a_3$ 
(see \cite[Pg.\ 406]{gkz94} and Remark \ref{rem:square}  
of Section \ref{sub:algor} below). Whether polynomial time complexity persists
for general $m$-nomials is already non-trivial for the tetranomial case since
the underlying discriminants become much more complicated. For example,
up to sign, $\Delta_{\{0,6,10,31\}}(c_1+c_2x^6+c_3x^{10}+c_4x^{31})$ is the
following homogeneous $21$-nomial of degree $35$:
\tiny
\[-2^{10} 3^{21} 5^{10} 7^{21} c^4_1c_3^{31}- 2^{13}
3^{24}5^{5}7^{21} c_1^2c_2^5c^{28}_3
-2^{14}3^{27}7^{21} c_2^{10}c^{25}_3+ 2^{8}3^{15}
5^{6}7^{16}\cdot 19\cdot 31^{6} 11731 c_1^7c_2^3 c_3^{23}c_4^2 \]
\[-2^{8}3^{13}\cdot 5\cdot  7^{13} \cdot 19\cdot 31^{4}\cdot 4931\cdot
11924843839 c_1^5c_2^8 c_3^{20}c_4^2 - 2^{6}3^{13}5^{6}7^{13}31^{13}
c_1^{12}c_2 c_3^{18}c_4^4 -2^{8}3^{16}5^{12}7^{13}\cdot 17\cdot 31^{2}
\cdot 3629537 c_1^3c_2^{13} c_3^{17}c_4^2 \]
\[-2^{6}3^{9}\cdot 5\cdot  7^{11}\cdot 29\cdot  31^{11} \cdot 6361\cdot
7477163 c_1^{10}c_2^6 c_3^{15}c_4^4 
-2^{10}3^{19}5^{22}7^{13} c_1c_2^{18} c_3^14c_4^2
-2^{6}3^{8}5^{12}7^{7}31^{9}\cdot 857\cdot 295269701 
c_1^8c_2^{11} c_3^{12}c_4^4 \]
\[ - 2^{7}3^8 \cdot 5\cdot  7^{7} 31^{18} \cdot 3327253
c_1^{15}c_2^4 c_3^{10}c_4^6 -2^{6}3 ^{4}5^{22}7^{4}31^7 160730667473
c_1^6c_2^{16} c_3^9c_4^4
+ 2^{4}3 ^{3}5^{12}7^{4}\cdot 11\cdot 29\cdot 31^{16}\cdot 233\cdot 1559
c_1^{13}c_2^9 c_3^7c_4^6 \]
\[-2^{6}5^{32}31^{5}\cdot 248552699041 c_1^4c_2^{21} c_3^6c_4^4
-2^{2}3 ^{4}57^4 \cdot 11\cdot 29\cdot  31^{25}
c_1^{20}c_2^2 c_3^5c_4^8 - 2^{5}\cdot 3 5^{22}31^{14} \cdot 113\cdot 317\cdot
461 c_1^{11}c_2^{14} c_3^4c_4^6 \]
\[- 2^{7}3 ^{3}5^{42}31^3 \cdot 41\cdot 367
c_1^2c_2^{26} c_3^3c_4^4
-2^{6}\cdot 3\cdot 5^{12}\cdot 19\cdot 31^{23}
c_1^{18}c_2^7 c_3^2c_4^8 + 2^4 3^4 5^{32} 31^{12} c_1^9c_2^{19} c_3c_4^6 -
31^{31} c^{25}_1 c_4^{10} - 2^{6}3^{6}5^{50} c^{31}_2 c^4_4 \]
                                                                                
\normalsize
\noindent
whose largest coefficient has $47$ digits. So we pose the following 
problem: 
\begin{disc} 
Given any finite subset $A\!\subset\!\Z$, define its 
{\bf size} to be the number of bits needed to write the binary expansions of 
all the points of $A$. Can one decide the vanishing of $A$-discriminants 
(for, say, specializations of the coefficients in $\C$) 
using a number of arithmetic operations polynomial in the 
size of $A$? \dia 
\end{disc} 
This highlights an embarassing gap in what is known about the complexity of 
discriminant computation: it is known that sparse
discriminants for {\bf bi}variate polynomials can be computed in polynomial
time on a Turing machine (resp.\ BSS machine over $\C$)
only if $\np\!\subseteq\!\pp$ (resp.\ $\np\!\subseteq\!\bpp$)
\cite[Main Theorem 1.4]{real}, but no such hardness result is known 
for the univariate case.
The preceding inclusions (see \cite{papa} for
a beautiful introduction to complexity theory) are currently considered quite
unlikely. On the positive side, it is known that 
$A$-discriminants (for arbitrary finite $A\!\subset\!\Zn$) can be decided 
in polynomial time when $A$ has $n+1$ or fewer points \cite[Prop.\ 1.8, pp.\ 
274--275]{gkz94}. Also, the vanishing of $A$-discriminants can be decided 
within the complexity class $\pp^{\np^\np}$ under a number-theoretic 
hypothesis strictly weaker than the Generalized Riemann Hypothesis \cite{dzh}. 

\section{Speed-Ups Through a Variant of $\alpha$-theory} 
\label{sec:back}
Checking whether a given point is an approximate root 
of a given polynomial can be done 
quite efficiently, thanks to the seminal work of Smale in \cite{smale}. 
Let us now formalize a refined version of this fact.  
\begin{dfn} 
For any analytic function $f : \C \longrightarrow \C$, let 
$\gamma(f,x)\!:=\!\sup\limits_{k\geq 2} \left| \frac{f^{(k)}(x)}
{k!f'(x)}\right|^{\frac{1}{k-1}}$. \dia 
\end{dfn} 
\begin{rem} 
\label{rem:gamma} 
It is worth noting that $1/\gamma(f,x_0)$ is a lower bound for the 
radius of convergence of the Taylor series of $f$ about $x_0$, so 
$\gamma(f,x_0)$ is finite whenever $f$ is nonsingular at $x_0$ 
\cite[Prop.\ 6, Pg.\ 167]{bcss}. \dia
\end{rem} 

\noindent
Recall that a function $g : U\longrightarrow \R$, defined on some
connected domain $U\!\subseteq\!\R$ is {\bf convex} iff
$g(\lambda x+(1-\lambda)y)\!\leq\!\lambda g(x)+(1-\lambda)g(y)$
for all $x,y\!\in\!U$ and $\lambda\!\in\![0,1]$. \dia
\begin{thm} 
\cite[Thm.\ 2]{ye}
\label{thm:ye} 
Suppose $f$ is convex and analytic on the open interval $(0,R)$. Also let 
$r\!\in\!(0,R)$ be any root of $f$ 
and suppose that there is an $\bar{\alpha}\!>\!0$ such that 
$x\gamma(f,x)\!\leq\!\bar{\alpha}$, 
for all $x$ in the closed interval 
$\left[\frac{r}{1+\frac{1}{8\bar{\alpha}}},\frac{r}{1-
\frac{1}{8\bar{\alpha}}}\right]$. 
Then for all $z\!\in\!(0,R)$: 
\begin{enumerate}
\item{If $f$ is monotonically decreasing in $[z,r]$ and $r\!\in\!\left[z,
(1+\frac{1}{8\bar{\alpha}})z\right]$ then $z$ is an \mbox{approximate} 
root of $f$.} 
\item{If $f$ is monotonically increasing in $[r,z]$ and $r\!\in\!\left[
(1-\frac{1}{8\bar{\alpha}})z,z\right]$ then $z$ is an \mbox{approximate} 
root of $f$. \qed } 
\end{enumerate} 
\end{thm} 
\begin{prop} 
\label{prop:sum} 
\cite[Prop.\ 1]{ye} 
Suppose $f_1 , f_2 : U \longrightarrow \R$ satisfy the hypotheses of  
Theorem \ref{thm:ye}, with $\bar{\alpha}_1$ and $\bar{\alpha}_2$ respectively 
playing the role of $\bar{\alpha}$, 
and that $f_1$ and $f_2$ are either both increasing or both decreasing. 
Then $f_1+f_2$ also satisfies the hypotheses of Theorem \ref{thm:ye}, 
with $\bar{\alpha}\!=\!\max\left\{\bar{\alpha}_1,\bar{\alpha}_2\right\}$. \qed 
\end{prop} 
\begin{prop} 
\label{prop:exp}
For any $D\!\in\!\R\!\setminus\!\{0\}$ and 
$c\!\in\!\R$, the function 
$f : [0,\infty) \longrightarrow \R$ defined by $x\mapsto x^D-c$ satisfies 
$x\gamma(f,x)\!=\!\left|\frac{D-1}{2}\right|$ for all $x\!>\!0$. \qed 
\end{prop} 

\noindent 
The second author observed the special case $D\!\in\!\Z\!\setminus\!\{0\}$ 
earlier in \cite[Example 1]{ye}. The proof their applies to aribtrary real $D$ 
as well. 

An additional technical result we will need is an analytic 
estimate which globalizes Theorem \ref{thm:ye}.  
But first let us observe an oscillation property we will need 
for our main algorithm. 
\begin{dfn} 
\label{dfn:global1} 
For any $D,m\!\in\!\N$ define $\cF(D,m)$ to be the family 
of $t$-nomials $f$ of degree $\leq\!D$ with $t\!\leq\!m$. For 
any $f\!\in\!\cF(D,m)\!\setminus\!\{0\}$, let $u_0(f)\!<\cdots<\!u_{N(f)}(f)$ 
be the ordered sequence comprised of $0$, the roots of $f'$ 
in $(0,\infty)$, and $\infty$. We then say that $f$ is {\bf dampened} 
iff for all $i\!\in\!\{0,\ldots,N(f)-1\}$, $f''$ has at most one 
root in $(u_i(f),u_{i+1}(f))$. \dia
\end{dfn} 

\noindent 
We of course have $\cF(D_1,m_1)\!\subseteq\!\cF(D_2,m_2)$ whenever 
$D_1\!\leq\!D_2$ and $m_1\!\leq\!m_2$. 
Rolle's Theorem applied to $f'$ tells us that $(u_i(f),u_{i+1}(f))$ 
(for $i\!\in\!\{1,\ldots,N(f)-2\}$) always contains at least one root of 
$f''$. It is then easy to show via Descartes' Rule and a routine 
calculation that $m$-nomials are always dampened when 
$m\!\leq\!4$. 

The invariant defined below generalizes the local quantity 
maximized in Theorem \ref{thm:ye} above, and thus helps enforce the 
accelerated convergence of Newton's method in certain cases of interest. 
\begin{dfn}
\label{dfn:global} 
For any $D,m\!\geq\!2$, any dampened $f\!\in\!\cF(D,m)$, and 
$i\!\in\!\{1,\ldots,N(f)-2\}$, define $v_i(f)$ to be the unique 
root of $f''$ in $(u_i(f),u_{i+1}(f))$, $v_{-1}(f)\!:=\!0$, and 
$v_{N(f)}(f)\!:=\!\infty$. Also define $v_0(f)$ to 
be $u_1(f)$ or the unique root of $f''$ in $(u_0(f),u_1(f))$,   
according as $f''$ is non-vanishing in $(u_0(f),u_1(f))$ or not.  
Also define $v_{N(f)-1}(f)$ to be $\infty$ or the unique root of 
$f''$ in $(u_{N(f)-1}(f),u_{N(f)}(f))$, according as $f''$ is non-vanishing 
in $(u_{N(f)-1}(f),u_{N(f)}(f))$ or not. Finally, define 
$\bar{\alpha}(D,m)$ to be 
\[ \sup\limits_{\substack{f \in \cF(D,m)\\ f \text{ dampened } \\ 
f(0)\neq 0 \text{ and } f'(0)= 0}} \ \ 
\sup\limits_{i\in\{0,\ldots,N(f)-1\}} 
\sup\limits_{\substack{x\in (v_{i-1}(f),v_i(f))\setminus\{u_i(f)\}\\ 
(v_{i-1}(f),v_i(f))\setminus\{u_i(f)\} \text{ contains a root of } 
f } } \left\{|x-u_i(f)|\gamma(f,x) \right\}. \text{ \dia} \] 
\end{dfn} 
While $\bar{\alpha}(D,m)$ may be hard to compute exactly, we can 
at least bound it above and below explicitly in some cases of interest. 
\begin{thm} 
\label{thm:global} 
For any $D,m\!\geq\!2$ we have 
$\bar{\alpha}(D,m)\!\leq\!\bar{\alpha}(D,m+1)$. 
In particular, \[\frac{D-1}{2}\!=\!\bar{\alpha}(D,2)\!\leq\!
\bar{\alpha}(D,3)\!\leq\!(D-1)(D-2)/2.\]  
\end{thm} 

One application of our $\bar{\alpha}$-invariant is the following generalization 
of Theorem \ref{thm:uni}. 
\begin{dfn} 
\label{dfn:trick} 
For any $m$-nomial $f$, let $\delta(f)$ be the smallest 
exponent appearing in $f$. Let us then define the operators 
$S : \cF(D,m) \longrightarrow \cF(D,m)$, 
$L_1 : \cF(D,m) \longrightarrow \cF(D-1,m-1)$, 
and $L_2 : \cF(D,m) \longrightarrow \cF(D-2,m-1)$ by 
$S(f)\!:=\!x^{-\delta(f)}f$, $L_1(f)\!:=\!\frac{d}{dx}S(f)$, 
and $L_2(f)\!:=\!\frac{d^2}{dx^2}S(f)$. 
\end{dfn} 
\begin{thm}
\label{thm:big} 
Following the notation of Theorem \ref{thm:uni} and Definition 
\ref{dfn:global}, suppose $f\!\in\!\cF(D,m)$ and that all the polynomials 
$\left\{(S\circ L_{e_1}\circ \cdots \circ L_{e_k})(f) \; | \; 
k\!\in\!\{0,\ldots,m-1\} , \ e_i\!\in\!\{1,2\}\right\}$ are dampened.  
Then one can find 
$\eps$-approximations to all the roots of $f$ in $[0,R]$ using just  
\[\cO\!\left(2^m\left\{\log(D)\log\left(\bar{\alpha}(D,m)\log 
\frac{R}{\eps}\right)+m^2 \cK(D,m)\right\}\right)\] 
arithmetic operations, where $\cK(D,m)$ is any upper bound on the worst-case 
arithmetic complexity of counting the roots of an arbitrary $f\!\in\!\cF(D,m)$ 
in an arbitrary interval. In particular, if $f\!\in\!\cF(D,4)$ then 
the polynomials\\
\mbox{}\hfill $\left\{S(f),S(L_{e_1}(f)),S(L_{e_1}(L_{e_2}(f))),S(L_{e_1}(L_{e_2}(L_{e_3}(f)))) \; | \; 
e_i\!\in\!\{1,2\}\right\}$\hfill\mbox{}\\  
are all dampened, for any $D\!\in\!\N$. 
\end{thm} 
It is easily checked (recalling our observations after Definition 
\ref{dfn:global1}) that any $m$-nomial with $m\!\leq\!4$ satisfies 
the hypothesis of our complexity bound above. 
Theorem \ref{thm:big} thus opens up the possibility of super-fast $m$-nomial 
solving when $m\!>\!3$ is fixed. In particular, good upper bounds 
on $\cK(D,m)$ and $\bar{\alpha}(D,m)$ appear to be a fundamental 
first step. 
\begin{probb} 
Is there an absolute contant $\kappa'$ such that 
$\bar{\alpha}(D,m)\!=\!O\!\left(\log^{\kappa'm} D\right)$? \dia  
\end{probb} 

\noindent 
An important related problem, in the spirit of \cite{mr04}, is how 
to bound from below the probability that a random $m$-nomial 
of degree $D$ (with $m\!\geq\!4$ fixed) is dampened. 

We state our underlying algorithms in the next section. We then prove
Theorems \ref{thm:uni}, \ref{thm:big}, and \ref{thm:bisturm} 
in Section \ref{sec:corcom}. We conclude with the proofs of
Theorems \ref{thm:tri} and \ref{thm:global} in Sections \ref{sub:tri}
and \ref{sub:hard}, respectively.

\section{The Algorithm and Subroutines}
\label{sub:algor} 
The central algorithm we use to prove Theorems \ref{thm:uni} and \ref{thm:big}, 
{\tt MNOMIALSOLVE}, is detailed below. 
Succinctly, the key idea behind {\tt MNOMIALSOLVE} is to subdivide the input 
interval into sub-intervals on which $\pm f$ is convex and monotonic, along 
with some additional sub-intervals on which $f$ is less well-behaved. 

\noindent
\mbox{}\hfill\epsfig{file=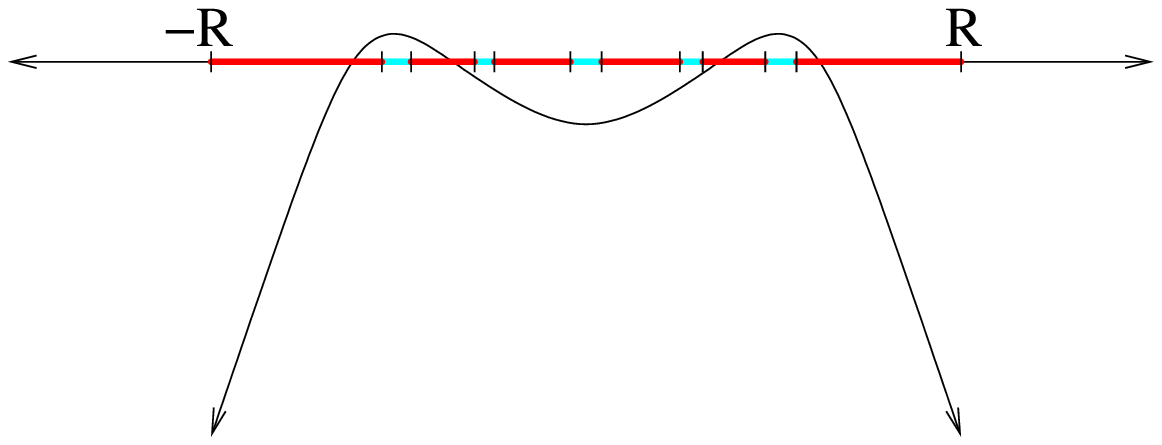,height=1.3in}\hfill\mbox{}\\
{\small {\sc Figure 1} $[-R,R]$ can be expressed as a union of $6$ intervals
on which $\pm f$ is convex and monotonic, and $5$ additional
intervals of width $<\!2\eps$ containing the roots of $f'f''$. \dia }

\medskip

\noindent 
Then, using 
special properties of sparse polynomials, a suitable combination of bisection 
and Newton iteration (detailed in the subroutine {\tt HYBRID}) yields 
$\eps$-approximations to all the roots in $(0,R)$. (Checking whether 
$0$ and/or $R$ are roots can be done simply by evaluating there.) Each 
approximation is 
guaranteed to correspond to its own unique root, and the roots in intervals on 
which $f$ behaves badly are certified via a subroutine called 
{\tt FASTERCOUNT}. The aforementioned subroutines are 
described shortly after our main algorithm.  

\medskip 

\small 
\vbox{
\noindent 
{\bf \fbox{ALGORITHM {\tt MNOMIALSOLVE}}}
\begin{itemize}
\item[{\bf INPUT}]{Real numbers $R$ and $\eps$ with
$0\!<\!\eps\!<\!R$, a univariate polynomial 
$f(x)\!:=\!c_1x^{a_1}+\cdots+c_{m}x^{a_{m}}\!\in\!\R[x]$ 
with exactly $m$ monomial terms (where 
$0\!\leq\!a_1\!<\cdots<\!a_m\!=\!D$) and (following the notation of 
Definition \ref{dfn:trick}) all the polynomials 
$\left\{(S\circ L_{e_1}\circ \cdots \circ L_{e_k})(f) \; | \; 
k\!\in\!\{0,\ldots,m-1\} , \ e_i\!\in\!\{1,2\}\right\}$ dampened, and 
an upper bound $\bar{\alpha}^*$ on $\bar{\alpha}(D,m)$.}  
\item[{\bf OUTPUT}]{A (possibly empty) multiset $Z\!\subseteq\!\R$, 
such that every root of $f$ in $[0,R]$ is 
$\eps$-approximated by a unique $z_i\!\in\!Z$, and $\#Z$ is exactly the number 
of roots of $f$ in $[0,R]$.} 
\end{itemize}} 

\noindent
{\bf DESCRIPTION}
\begin{itemize} 
\item[{\bf Step 0}]{Set $Z\!:=\!\emptyset$.} 
\item[{\bf Step 1}]{If $m\!=\!0$ then set $Z\!:=\!\R$ and STOP.} 
\item[{\bf Step 2}]{If $m\!=\!1$ and $a_1\!>\!0$ then set $Z\!:=\!\{0\}$ and 
STOP.} 
\item[{\bf Step 3}]{If $m\!=\!1$ and $a_1\!=\!0$ then STOP.} 
\item[{\bf Step 4}]{If $a_1\!>\!0$ then set 
$Z\!:=\!\{0\}$, $f\!:=\!f/x^{a_1}$, $a_i\!:=\!a_i-a_1$ for 
all $i\!\geq\!2$, and then $a_1\!:=\!0$.} 
\item[{\bf Step 5}]{If $a_m\!=\!1$ then set  
$Z\!:=\!\{-c_1/c_2\}$ or $\emptyset$, according as $c_1c_2\!\leq\!0$ 
or not, then STOP.}  
\item[{\bf Step 6}]{If $a_m\!=\!2$ then make the change of 
variables $x\!:=\!x-\frac{c_2}{2c_3}$.} 
\item[{\bf Step 7}]{If $f(R)\!=\!0$ then append $R$ to $Z$. } 
\item[{\bf Step 8}]{If $m\!=\!2$ and $f(0)f(R)\!\geq\!0$, then STOP.} 
\item[{\bf Step 9}]{If $m\!=\!2$ and $f(0)f(R)\!<\!0$, find 
an $\eps$-approximation $z$ to the unique root of $f$ in $(0,R)$ via 
subroutine {\tt HYBRID} (with input $(\eps,R,\pm f,\bar{\alpha}^*)$, where 
the sign is chosen so that $f$ is convex on $(0,R)$), set 
$Z\!:=\!Z\cup \{z\}$, and STOP.} 
\item[{\bf Step 10}]{Otherwise, if $a_2\!>\!1$, find via algorithm {\tt 
MNOMIALSOLVE} an ordered sequence $u_0\!<\cdots<\!u_{k_1}$ comprised 
of $0$, a set of $\eps$-approximations (of the correct 
cardinality) in $(0,R)$ for the roots of $f'$ in $(0,R)$, and $R$. 
Then GOTO Step 12.} 
\item[{\bf Step 11}]{If $a_2\!=\!1$ then redefine $f(x)\!:=\!x^{a_3} f(1/x)$ 
and $\eps\!:=\!\min\left\{1,\frac{\eps^3}{2},\frac{1}{2R^2}\right\}$ and 
make the change of variables $x\!:=\!1/x$. } 
\item[{\bf Step 12}]{Define $v_{k_1}\!:=\!R$.} 
\item[{\bf Step 13}]{Using algorithm {\tt FASTERCOUNT}, 
decide if $f''$ has a root in $[u_0+\eps,u_1-\eps]$. If so, 
let $L\!:=\!u_1-u_0-2\eps$. Then, evaluating $f''$ at 
$u_0+\eps$ and $u_1-\eps$, and using algorithm {\tt HYBRID} with 
input $(\eps,L,\pm f,\alpha^*)$ (the sign chosen so that $\pm f$ is 
convex on $(u_0+\eps,u_1-\eps)$), define $v_0$ to be 
the $\eps$-approximation found for the unique root of $f''$ in 
$[u_0+\eps,u_1-\eps]$. Otherwise, define $v_0$ to 
be $u_1$. } 
\item[{\bf Step 14}]{ Using algorithm {\tt FASTERCOUNT},
decide if $f''$ has a root in $[u_{k_1-1}+\eps,u_{k_1}-\eps]$. If so,
let $L\!:=\!u_{k_1}-u_{k_1-1}-2\eps$, $g(x)\!:=\!f(x+u_{k_1-1})$. 
Then, evaluating $f''$ at $u_{k_1-1}+\eps$ and $u_{k_1}-\eps$, and 
using algorithm {\tt HYBRID} with
input $(\eps,L,\pm g,\alpha^*)$ (the sign chosen so that $\pm g$ is
convex on $(u_{k_1-1}+\eps,u_{k_1}-\eps)$), define $v_{k_1-1}$ to be
the $\eps$-approximation found for the unique root of $f''$ in
$[u_{k_1-1}+\eps,u_{k_1}-\eps]$. Otherwise, define $v_{k_1-1}$ to
be $R$.}  
\item[{\bf Step 15}]{For all $i\!\in\!\{1,\ldots,k_1-1\}$, count 
via algorithm {\tt FASTERCOUNT} the number $m_i$ 
of roots of $f$ in $(u_i-\eps,u_i+\eps) \cap (0,R)$ closest to $u_i$ 
and no other root of $f'f''$, and append $m_i$ copies of $u_i$ to $Z$.}  
\item[{\bf Step 16}]{For all $i\!\in\!\{0,\ldots,k_1-1\}$, count via algorithm 
{\tt FASTERCOUNT} the number $n_i$ 
of roots of $f$ in $(v_i-\eps,v_i+\eps) \cap (0,R)$ closest to $v_i$ and no 
other root of $f'f''$, and append $n_i$ copies of $v_i$ to $Z$.}  
\item[{\bf Step 17}]{FOR $i\!\in\!\{0,\ldots,k_1-1\}$ DO  
\begin{itemize}
\item[{\bf Step 17(a)}]{Evaluate $f$ at $u_i+\eps$ and  
$v_i-\eps$. If $f(u_i+\eps)\!=\!0$ (resp.\ 
$f(v_i-\eps)\!=\!0$) then append $u_i+\eps$ (resp.\ $v_i-\eps$) to 
$Z$ and GOTO Step 17(c). } 
\item[{\bf Step 17(b)}]{If $f(u_i+\eps)f(v_i-\eps)\!<\!0$ then let 
$L_i\!:=\!v_i-u_i-2\eps$ and $g_i(x)\!:=f(x+u_i+\eps)$. Then, using 
subroutine {\tt HYBRID} with input $(\eps,L_i,\pm g_i,\bar{\alpha}^*)$ 
(the sign chosen so that $g_i$ is convex on $(0,v_i-u_i-2\eps)$), 
find an $\eps$-approximation in $(u_i+\eps,v_i-\eps)$ to 
the unique root of $f$ in $(u_i+\eps,v_i-\eps)$, and append this 
approximation to $Z$.} 
\item[{\bf Step 17(c)}]{Evaluate $f$ at $v_i+\eps$ and
$u_{i+1}-\eps$. If $f(v_i+\eps)\!=\!0$ (resp.\
$f(u_{i+1}-\eps)\!=\!0$) then append $v_i+\eps$ (resp.\ $u_{i+1}-\eps$) to
$Z$ and GOTO Step 18. } 
\item[{\bf Step 17(d)}]{If $f(v_i+\eps)f(u_{i+1}-\eps)\!<\!0$ then let
$L_i\!:=\!u_{i+1}-v_i-2\eps$ and $g_i(x)\!:=f(u_{i+1}-\eps-x)$.
Then, using subroutine {\tt HYBRID} with input 
$(\eps,L_i,\pm g_i,\bar{\alpha}^*)$ 
(the sign chosen so that $g_i$ is convex on $(0,u_{i+1}-v_{i}-2\eps)$), 
find an $\eps$-approximation in $(v_{i}+\eps,u_{i+1}-\eps)$ to
the unique root of $f$ in $(v_{i}+\eps,u_{i+1}-\eps)$, and append this 
approximation to $Z$. }  
\end{itemize}} 
\item[{\bf Step 18}]{END FOR} 
\item[{\bf Step 19}]{STOP} 
\end{itemize} 
\normalsize 

\begin{rem} 
\label{rem:k}
In reality, one should replace the use of subroutine {\tt FASTERCOUNT} 
above by the best algorithm attaining the complexity bound of $\cK(D,m)$ 
alluded to in Theorem \ref{thm:big}, whenever $m\!\geq\!4$. 
We have omitted this detail above simply to keep the algorithm definite. 
\dia  
\end{rem} 

\medskip 

\scriptsize 
\noindent 
{\bf \fbox{SUBROUTINE {\tt HYBRID} (Compare \cite[Pg.\ 277]{ye})}}
\begin{itemize} 
\item[{\bf INPUT}]{$\eps,R\!\in\!\R$ (with $0\!<\!\eps\!<\!R$), 
a monotonic analytic function $\phi : (0,R)\longrightarrow \R$ (with 
$\phi(\eps)\phi(R)\!<\!0$ and $\phi$ convex on $(0,R)$), and a  
positive upper bound $\bar{\alpha}^*$ on $z\gamma(\phi,z)$ valid for 
all $z\!\in\!(0,R)$ . } 
\item[{\bf OUTPUT}]{An $\eps$-approximation of the unique root of $\phi$ in 
$(0,R)$, using no more than 
$O\!\left(\log\left(\alpha\log\frac{R}{\eps}
\right)\right)$ evaluations of $\phi$ and $\phi'$, and 
$O\!\left(\log\left(\alpha\log\frac{R}{\eps}
\right)\right)$ additional arithmetic operations. }
\end{itemize} 

\noindent
{\bf DESCRIPTION} 
\begin{itemize} 
\item[{\bf Step 0}]{Define $c_0$ to be $1+\frac{1}{8\balpha}$ 
or $\frac{1}{1-\frac{1}{8\balpha}}$, according as $\phi$ is 
decreasing or increasing on $(0,R)$. Then define $c_{i}\!:=\!c^2_{i-1}$ for 
all $i\!>\!0$, and compute $c_1,\ldots,c_M$ where $M$ is the first positive 
integer with $c_M\!=\!c^{2^M}_0\!\geq\!\frac{R}{\eps}$. Finally, set 
$\hat{x}\!:=\!\eps$ and $\hat{k}\!:=\!M$.} 
\item[{\bf Step 1}]{ If [($\phi$ is decreasing on $(0,R)$ and 
$\phi(c_{\hat{k}-1}\hat{x})\!>\!0$) or 
($\phi$ is increasing on $(0,R)$ and
$\phi(c_{\hat{k}-1}\hat{x})\!<\!0$)] then 
set $\hat{x}\!:=\!c_{\hat{k}-1}\hat{x}$ and $\hat{k}\!:=\!\hat{k}-1$ 
and GOTO Step 1. } 
\item[{\bf Step 2}]{ Otherwise, if $k\!>\!0$ then set 
$\hat{k}\!:=\!\hat{k}-1$ and GOTO Step 1.} 
\item[{\bf Step 3}]{Perform $\log_2\left(3+\log_2\frac{R}{\eps}\right)$ 
iterations of Newton's method (with $\hat{x}$ as the starting point), then 
OUTPUT the very last iterate.}  
\end{itemize} 

\normalsize 
\begin{rem}
\label{rem:square} 
A well-known trick we will use frequently and implicitly is the 
computation of $x^d$, where $x\!\in\!\R$ and $d\!\in\!\Z$, within 
$2\lceil\log_2(|d|+1)\rceil$ multiplications using 
$\leq\!\lceil\log_2(|d|+1)\rceil$ intermediate real numbers. Briefly, assuming 
the binary expansion of $d$ is $(a_k\cdots a_0)_2$ (so $a_0$ is the 1's bit), 
we simply compute $x^{2^i}$, via recursive squaring, for all $i$ with $a_i$ 
nonzero, and then multiply the resulting numbers together. For instance, 
$363\!=\!2^8 + 2^6 + 2^5 + 2^3 + 2 + 1\!=\!(101101011)_2$ and thus 
$x^{363}\!=\!(\cdots (x\underset{8}{\underbrace{^2)\cdots)^2}} 
\cdot (\cdots (x\underset{6}{\underbrace{^2)\cdots )^2}}\cdot 
(\cdots (x\underset{5}{\underbrace{^2)\cdots )^2}} \cdot ((x^2)^2)^2 
\cdot x^2 \cdot x$, which only requires $17$ multiplications 
and enough memory for $9$ intermediate real numbers --- much faster 
than the naive $362$ multiplications. \dia 
\end{rem} 
\scriptsize 

\medskip 

\vbox{
\noindent 
{\bf \fbox{SUBROUTINE {\tt FASTERCOUNT}}} 
\begin{itemize}
\item[{\bf INPUT}]{A polynomial $p\!\in\!\R[x_1]$ with exactly $m\!\geq\!3$  
monomial terms, and $a,b\!\in\!\R$ with $a\!<\!b$. }
\item[{\bf OUTPUT}]{The number of roots of $p$ in the open interval 
$(a,b)$.}
\end{itemize} 

\noindent
{\bf DESCRIPTION} 
\begin{itemize}
\item[{\bf Step 0}]{If $m\!>\!3$ then use 
Sturm-Habicht sequences (as in, say, \cite{lickroy}) to 
count the number of roots in $(a,b)$ and STOP.}   
\item[{\bf Step 1}]{Otherwise, let $p_0\!:=\!p$, $p_1\!:=\!p'$, 
let $-p_2$ be the remainder of $p_0/p_1$, and set $i\!:=\!2$. } 
\item[{\bf Step 2}]{If $p_i\!\in\!\R$ then set $K\!:=\!i$ and 
GOTO Step 5. } 
\item[{\bf Step 3}]{Write $p_i(x)\!:=\!u_1 x^{a_1}
-u_0x^{a_0}$, $p_{i-1}(x)\!:=\!v_1 x^{b_1}-v_0 x^{b_0}$ 
(where $a_1\!>\!a_0$ and $b_1\!>\!b_0$), and let $c_j\!:=\!\left\lceil 
\frac{b_j-a_1}{a_j-a_0}\right\rceil$ for $j\!\in\!\{0,1\}$. } 
\item[{\bf Step 4}]{ Replace $i$ by $i+1$, define $p_i(x)\!:=\!-v_1
(\frac{u_1}{u_0})^{c_1}x^{b_1-c_1(a_1-a_0)}+v_0 
(\frac{u_1}{u_0})^{c_0}x^{b_0-c_0(a_1-a_0)}$, and GOTO Step 2.} 
\item[{\bf Step 5}]{ Using recursive squaring (cf.\ Remark \ref{rem:square}), 
evaluate and RETURN $N_A-N_B$ (cf.\ Definition \ref{dfn:sturm}) 
where $A$ (resp.\ $B$) is $(p_0(a),\ldots,p_K(a))$ (resp.\ 
$(p_0(b),\ldots,p_K(b))$). } 
\end{itemize}  } 
\normalsize 

\section{Correctness and Complexity: Proving Theorems \ref{thm:uni}, 
\ref{thm:big}, and \ref{thm:bisturm}} 
\label{sec:corcom} 

\noindent
{\bf Proof of Theorem \ref{thm:uni}:}\\  
For open intervals, the second assertion follows immediately from Lemma \ref{lemma:sturm} and Theorem \ref{thm:bisturm}. The case of closed intervals 
follows almost identically, save for an additional evaluation of $f$ at the 
end-points of the interval which increases the complexity bound by (a 
negligible) $O(\log D)$. 

The first assertion follows immediately from 
Theorems \ref{thm:big} and \ref{thm:global}, and the second assertion 
of Theorem \ref{thm:uni} which we've just proved. \qed 

\medskip
\noindent
{\bf Proof of Theorem \ref{thm:big}:}\\ 
Henceforth, for convenience, we will say ``time'' in place of 
``arithmetic complexity.'' 
Clearly, it suffices to show that 
{\tt MNOMIALSOLVE} is correct and satisfies the stated complexity bound 
when we set $\balpha\!=\!\bar{\alpha}(D,m)$.  
Toward this end, let us first observe that 
$L_1(f)$ has exactly $m-1$ monomial terms and $L_2(f)$ has
$\leq\!m-1$ monomial terms. Descartes' Rule then implies that $L_1(f)$
and $L_2(f)$ each have no more than $m-2$ positive roots.
Note also that $f$, $f'$, and $f''$ can all be evaluated using just
$O(m\log D)$ arithmetic operations, thanks to Remark \ref{rem:square}. 
These observations will be used implicitly throughout our proof. 

The correctness of {\tt MNOMIALSOLVE} is then straightforward for 
$m\!\leq\!2$. Also, for $m\!=\!2$ and $D\!\geq\!3$, the complexity of 
{\tt MNOMIALSOLVE} is clearly the same (asymptotically) as that of subroutine 
{\tt HYBRID} with $\balpha\!:=\!\frac{D-1}{2}$, thanks to Proposition 
\ref{prop:exp}. So we can assume $m\!\geq\!3$. 

In the special case where $D\!=\!2$ and $m\!=\!3$, note that Step 6 is 
but an implementation of completing the square. So, provided one 
shifts one's $\eps$-approximations by $\frac{c_2}{2c_3}$ while  
running the remainder of {\tt MNOMIALSOLVE}, we can clearly attain the stated 
complexity bound for $D\!=\!2$. So we can assume additionally that 
$D\!\geq\!3$ and now focus on Steps 10--18. 

In the special case $a_2\!=\!1$, note that the change 
of variables $x\mapsto 1/x$ maps the interval $(r-\delta,r+\delta)$ 
to $\frac{1}{r^2-\delta^2}(r-\delta,r+\delta)$ for 
$r\!>\!\delta\!>\!0$. An elementary 
calculation then reveals that $\delta\!\leq\!\max\left\{1,\eps^3/2\right\} 
\Longrightarrow$ for any $\delta$-approximation $y$ of $x\!\geq\!\eps$,  
$1/y$ is an $\eps$-approximation of $1/x$. Note also that 
$x^{a_3}f(1/x)\!=\!c_3+c_2x^{a_3-a_2}+c_1x^{a_3}$ (since $a_1\!=\!0$ 
after Step 4), and $a_3-a_2\!>\!1$ since $a_2\!=\!1$ and $a_3\!=\!D\!\geq\!3$. 
So, provided one takes a reciprocal after running {\tt MNOMIALSOLVE}, we see 
that the worst-case time complexity of any instance with 
$a_2\!=\!1$ is asymptotically no worse than that of the instances with 
$a_2\!>\!1$. So we can assume additionally that $a_2\!>\!1$.

Now, by assumption, $S(f)$ is dampened, as are all the inputs 
going into recursive calls of {\tt MNOMIALSOLVE}. Note also that 
by our dampening assumption (and Rolle's Theorem), the intervals 
from Steps 17(b) and 17(d) indeed contain exactly one root of $f''$. 
So our algorithm is well-defined and indeed finds $\eps$-approximations 
to all the roots of $f$ in $[0,R]$, assuming 
subroutines {\tt FASTERCOUNT} and {\tt HYBRID} are correct. 
The correctness of the latter two subroutines is proved below, so let us now 
concentrate on proving our main complexity bound. 

So let $\cC(m)$ denote the time needed for {\tt MNOMIALSOLVE} 
to execute completely for an input consisting of an $f\!\in\!\cF(D,m)$ 
(satisfying the dampening assumptions of our current theorem), an interval 
length of $\leq\!R$, and a precision $\eps$. Clearly, 
assuming subroutine {\tt HYBRID} runs in time 
$\cO\!\left(m\log(D)\log\left(\bar{\alpha}^*\log 
\frac{R}{\eps}\right)\right)$ (which is covered in a separate proof 
below), we then have 
that $\cC(m)$ must satisfy the following recurrence relation: 
\[\cC(m)\!\leq\!2\cC(m-1)
+ \cO\!\left(m(m\log D)\log\left(\bar{\alpha}^*\log\frac{R}{\eps}
\right)\right)+(2m-4)\cK(D,m),\]  
where the first term corresponds to finding the 
$\eps$-approximations of the roots of $f'$ and $f''$ (Steps 10, 13, and 
14),  the last term corresponds to the application of algorithm {\tt 
FASTERCOUNT} to the intervals about these roots (Steps 13--16), 
and the $\cO$ term corresponds to the application of 
algorithm {\tt HYBRID} to the intervals between these roots (Steps 17--18). 
(Note that we have implicitly used the fact that $\cK(D,m)$ and 
$\alpha(D,m)$ are non-decreasing functions of $D$ and $m$.)  
In particular, for $m\!\geq\!4$, the last term of our recurrence 
is justified by Remark \ref{rem:k}. 

Noting that we can regard $R$ and $\eps$ as constants, it is then 
clear that to find the asymptotics of our recurrence, it 
suffices to find the asymptotics of the recurrence 
\[c_{m}\!\leq\!2c_{m-1} +m^2A+(2m-4)B(m).\] 
A simple {\tt Maple} calculation 
\mbox{({\tt rsolve(c(m)=2*c(m-1)+a*m*(m-1)+b(m)*(2*m-4),c(m));})}  
then shows that any such $c_m$ must satisfy 
$c_m\!=\!O(2^m(c_2+A)+\sum^m_{j=1}
2^{m-j}jB(j))$. So, knowing that $\cC(2)\!=\!O\!\left(\log(D)\log 
\log\frac{R}{\eps}\right)$ \cite[Pg.\ 279]{ye}, we must then have that 
\[ \cC(m)=O\!\left(2^m\left\{\log(D)\log\left(\bar{\alpha}^*\log \frac{R}{\eps}
\right)+m^2\cK(D,m)\right\}\right). \text{ \qed} \] 

\noindent 
{\bf Proof of Theorem \ref{thm:bisturm}:}\\ 
The complexity bound follows immediately from the upper bound on $K$, 
assuming algorithm {\tt FASTERCOUNT} is correct. The bound on $K$ and the 
correctness of algorithm {\tt FASTERCOUNT} are proved below. \qed 

\medskip 
 
At this point, we need only show that the subroutines {\tt HYBRID} and 
{\tt FASTERCOUNT} are correct, and that {\tt FASTERCOUNT} satisfies a 
sufficiently good complexity bound.  

\medskip 
\noindent 
{\bf Proof of Correctness and Complexity Analysis of {\tt HYBRID}:}\\ 
First note that by assumption, there is a unique $j\!\in\!\{0,\ldots,2^{K-1}\}$ 
such that $(\beta^j\eps,\beta^{j+1}\eps]$ contains the sole root of 
$\phi$ in $(0,R)$. It then becomes clear that Steps 0--2 are merely an 
implementation of bisection that finds this $j$: In particular, 
Steps 0--2 essentially find the binary expansion of $j$ (from 
most significant bit to least significant bit), and $\hat{x}\!=\!
\beta^j$ once Step 3 is reached. By Theorem \ref{thm:ye}, $\hat{x}$ is then 
an approximate root of $\phi$, and Step 3 indeed finds an 
$\eps$-approximation of our root. 

As for the complexity of our algorithm, it is clear that 
$M\!=\!\left\lceil\log_2
\left(\log_2\frac{R}{\eps}\right)-\log_2\log_2\beta\right\rceil$ and 
we thus need only 
$O\!\left(\log\left(\log\frac{R}{\eps}\right)-\log\log\beta\right)$ 
evaluations of $\phi$ (and inequality checks) until Step 3. Step 3 then 
clearly takes just $O\!\left(\log \log \frac{R}{\eps}\right)$ arithmetic 
operations and evaluations of $\phi$ and $\phi'$. 

To conclude, recall the elementary inequalities 
$\frac{\pm 1}{\log(1\pm x)}\!\leq\!x$, valid for all $x\!\geq\!0$. 
Observe then that 
$-\log\log \beta\!=\!\log \frac{1}{\log\beta}\!\leq\!
\log \balpha$. Since $\phi$ an $m$-nomial of degree $D$ implies that 
$\phi$ and $\phi'$ can be evaluated within $O(m\log D)$ arithmetic 
operations (cf.\ Remark \ref{rem:square}), we are done. \qed 

\medskip
\noindent
{\bf Proof of Correctness and Complexity Analysis of {\tt FASTERCOUNT}:}\\ 
The correctness of {\tt FASTERCOUNT} when $m\!>\!3$ is clear from 
any of the standard references \cite{marie,lickroy}  
so there is nothing to prove. 

For $m\!=\!3$, correctness 
follows easily from Lemma \ref{lemma:sturm} upon noting a simple 
fact: The remainder of $\frac{v_1x^{b_1}-v_0x^{b_0}}
{u_1x^{a_1}-u_0x^{a_0}}$ is nothing more than the reduction 
of $v_1x^{b_1}-v_0x^{b_0}$ modulo $u_1x^{a_1}-u_0x^{a_0}$. 
In particular, Step 2 of {\tt FASTERCOUNT} is just a compact 
representation of a sufficient number of applications of the identity 
$x^{a_1}=\frac{u_0}{u_1}x^{a_0}$ to reduce 
$-(v_1x^{b_1}-v_0x^{b_0})$ modulo $u_1x^{a_1}-u_0x^{a_0}$ to 
a polynomial of degree $<\!a_1$. So correctness follows immediately. 

As for the complexity bound, we need only observe that 
(by recursive squaring) every execution of Step 2 and Step 3 takes only 
$\cO(\log D)$ arithmetic operations, and that $K\!=\!\cO(\log D)$ 
at the termination of the algorithm, i.e., there are only $\cO(\log D)$ 
remainders in our Sturm sequence. The first assertion is clear, so to prove 
the latter assertion we will need to prove a technical bound on the 
exponents which occur in our remainder sequence. 

In particular, let $\ell_i$ be the absolute value of the difference exponents 
of $p_i$ (we set $\ell_i\!:=\!0$ if $p_i$ is monomial). Note that  
$p_i$ a monomial $\Longrightarrow p_{i+1}$ is a monomial, and thus 
$p_{i+2}$ is constant. Note also that, via the definition of long division, 
we have $\ell_{i+2}\!\leq\!|\ell_{i+1}-\ell_i|$, with equality occurring iff 
$p_{i+2}$ is a binomial. So if we can show 
\[ (\star) \ \ \ \ \ \ \ \ \ \ \ \ \ \ \ \ \ \ \ \ \ \ \ \ \ \ \ \ \ \ \ \ \ 
\ \ \ \min\{\ell_{i+2},\ell_{i+3}\}\!\leq\!\ell_i/2\hspace{2in}\mbox{} \] 
for all $i$, then we will easily obtain our bound on $K$ 
(since $\ell_0,\ell_1\!\leq\!D$ and all the $\ell_i$ are integers). 

To prove ($\star$) observe that $\ell_{i+1}\!\leq\!\frac{1}{2}\ell_i 
\Longrightarrow \ell_{i+2}\!\geq\!\frac{1}{2}\ell_{i+1}$ (if $\ell_{i+2}\!\neq\!0$) 
and thus $\ell_{i+3}\!\leq\!\frac{1}{2}\ell_i$. Similarly, 
$\ell_{i+1}\!\geq\!\frac{1}{2}\ell_i \Longrightarrow \ell_{i+2}\!\leq\!\frac{1}{2}\ell_i$. 
So $K\!\leq\!3\lceil\log_2 D\rceil+2\!=\!\cO(\log D)$ and we are done. \qed 

\section{The Proof of Theorem \ref{thm:tri}} 
\label{sub:tri} 
Dividing by a suitable monomial term, we can clearly assume without 
loss of generality that $f$ has a constant term . 
Note then that $C$ is diffeomorphic to a line iff $f$ is not 
the square of binomial, via \cite[Prop.\ 2 and Lemma 1]{tri}. So there are 
actually no isolated singularities for trinomial curves.  

It is then clear that $x$ is a vertical tangent of $C \Longrightarrow 
x$ is a root of $F_v\!:=(f,\frac{\partial f}{\partial x_2})$. 
The roots of $F_v$ with $x$ or $y$ coordinate equal to $0$ or 
$R$ can be found easily by suitably fixing one of the variables 
and applying Theorem \ref{thm:uni}, so let us concentrate 
on the roots of $F_v$ in $(0,R)^2$. 

By \cite[Lemma 4.2]{real} and \cite{vanderk} we can then find (for 
any constant $\delta\!>\!0$) a monomial change 
of variables, using just $\cO(\log^{2+\delta} D)$ bit operations 
and $\cO\!\left(\log(D)\log\left(D\log\frac{R}{\eps}\right)\right)$ arithmetic 
operations, such that the resulting sparse polynomial system $G\!:=\!(g_1,g_2)$ 
has $g_1\!\in\!\R[x_1,x_2]$ a trinomial, $g_2\!\in\!\R[x_1]$ a binomial, and 
all underlying exponents have bit-length $\cO(\log^2 D)$. We can then solve 
$g_2$ to accuracy $\min\left\{\eps^{\cO(D^2)},1\right\}$ 
via Theorem \ref{thm:uni}, back-substitute the 
roots into $g_1$, solve the univariate specialized $g_1$ to accuracy 
$\min\left\{\eps^{\cO(D^2)},1\right\}$, then invert our monomial map via 
\cite[Main Theorem 1.3]{real} and \cite{vanderk} to 
recover the roots of $F_v$, using just $O(\log^{2+\delta}D)$ additional bit 
operations. So we can indeed find isolated vertical tangents 
within the stated complexity bound.  

To conclude, note that \cite[Lemma 9]{tri} tells us that 
if $x$ is an isolated inflection point of $C$ then 
$[\partial^2_1 f\cdot  (\partial_2 f)^2-
2\partial_1\partial_2 f\cdot\partial_1f\cdot\partial_2f
+\partial^2_2 f\cdot(\partial_1 f)^2]\!=\!0$, where 
$\partial_i\!:=\!\frac{\partial}{\partial x_i}$. In the special case of a 
trinomial, say $f(x_1,x_2)\!:=1+Ax^a_1x^b_2+Bx^c_1x^d_2$, the 
preceding polynomial in derivatives is exactly: 
\[-ab\left (a+b\right ){R}^{3}+\left ({b}^{2}{c}^{2}+{a}^{2}{d}^{2}-b^2 c-
2\,abd-{a}^{2}d-2\,abc-2\,abcd\right ){R}^{2}S\]
\[+\left (b^2 c^2-2\,abcd-a{d}^{2}-2\,acd-2\,bcd-b{c}^{2}+{a}^{2}{d}^{2}
\right)RS^2-cd\left (c+d\right )S^3,\]
\noindent
where $R\!:=\!Ax^a_1x^b_2$ and $S\!:=\!Bx^c_1x^d_2$,  
thanks to a {\tt Maple} calculation. (Note that we must also have 
$1+R+S\!=\!0$ 
since we're assuming $f\!=\!0$.) So we efficiently solve for $x$ as follows: 
(1) solve the above cubic for $\frac{R}{S}$ up to 
accuracy $\min\left\{\eps^{\cO(D)},1\right\}$ via the main algorithm of 
\cite{neffreif} or \cite{binipan}, (2) solve for $(R,S)$ via the 
additional relation $1+R+S\!=\!0$, (3) solve the resulting binomial 
systems in $(R,S,x_1,x_2)$ for $x$ via \cite[Main Theorem 1.3]{real} 
and \cite{vanderk}. So we can $\eps$-approximate the inflection points of $C$ 
within the stated asymptotic complexity bound as well, and we are done. \qed  

\section{The Proof of Theorem \ref{thm:global}}
\label{sub:hard} 
\begin{dfn}
For any $a\!\in\!\R$ and $k\!\in\!\N\cup\{0\}$, let 
$(a)_k\!:=\!a(a-1)\cdots (a-k+1)$. (So  $(a)_0\!=\!1$, 
$(a)_1\!=\!a$, $a\!\in\!\N \Longrightarrow 
(a)_a\!=\!a!$, and $a,k\!\in\!\N$ with $k\!>\!a \Longrightarrow 
(a)_k\!=\!0$.)  
\dia 
\end{dfn}   
\noindent 

\noindent 
{\bf Proof of Theorem \ref{thm:global}:}\\ 
The first bound is immediate since $\cF(D,m)\!\subseteq\!\cF(D,m+1)$.

The formula for $\bar{\alpha}(D,2)$ follows easily since $m\!=\!2$ and 
$f(0)\!\neq\!0$ implies that $f'$ has no positive roots. 
This in turn implies that $N(f)\!=\!1$, so the only interval we (may)  
need to majorize over is $(u_0(f),u_1(f))\!=\!(0,\infty)$. In particular,  
the quantity we majorize is just $x\gamma(f,x)$, and the formula 
for $\bar{\alpha}(D,2)$ then follows immediately from Proposition 
\ref{prop:exp}. So we can assume $m\!=\!3$.  
Since $\gamma(f,x)\!=\!\gamma(cf,x)$ for any nonzero constant $c$, 
we can also clearly assume that $f(x)\!=\!x^{a_3}-Ax^{a_2}+B$ where 
$A$ and $B$ are nonzero real constants and $a_3\!=\!D\!>\!a_2\!\geq\!1$. 
Since we assume $f'(0)\!=\!0$ in the definition of $\bar{\alpha}(D,m)$, 
we must also have $D\!\geq\!3$ and $a_2\!>\!1$. 

Next, note that $A\!<\!0 \Longrightarrow f'$ has no positive roots. So, 
similar to the $m\!=\!2$ case, one either majorizes over the empty set 
or $(0,\infty)$. In the former case there is no contribution to 
$\bar{\alpha}(D,m)$, while in the second case, we 
see that we are again majorizing $x\gamma(f,x)$ over $(0,\infty)$. 
Since the latter supremum is no greater than $\frac{D-1}{2}$ by 
Propositions \ref{prop:sum} and \ref{prop:exp}, we can 
can now assume additionally that $A\!>\!0$. 

Now note that 
\begin{eqnarray*}
f'(x)& =& a_3x^{a_3-1}-a_2Ax^{a_2-1}\\
& = & a_3x^{a_2-1}\left(x^{a_3-a_2}-\frac{a_2A}{a_3}\right),
\end{eqnarray*}  
and
\begin{eqnarray*}
f''(x)& =& a_3(a_3-1)x^{a_3-2}-a_2A(a_2-1)x^{a_2-2}\\
& = & 
a_3(a_3-1)x^{a_2-2}\left(x^{a_3-a_2}-\frac{a_2A(a_2-1)}{a_3(a_3-1)}\right). 
\end{eqnarray*}  
Clearly then, 
$x_1 := \left(\frac{a_2A}{a_3}\right)^{1/(a_3-a_2)}$
is the unique positive root of $f'$ and 
$x_2 := \left(\frac{a_2(a_2-1)A} {a_3(a_3-1)}\right)^{1/(a_3-a_2)}$
is the unique positive root of $f''$. Furthermore, since 
$1\!<\!a_2\!<\!a_3$, it is clear that $0\!<\!x_2\!<\!x_1$. 
In particular, we see that we are left with a collection of $3$ intervals to 
consider in the definition of $\bar{\alpha}(D,m)$: 
\[ \left\{(0,x_2), (x_2,x_1), (x_1,\infty)\right\} \]  
and the quantities we majorize on these intervals are, respectively: 
\[ \left\{x\gamma(f,x),(x_1-x)\gamma(f,x),(x-x_1)\gamma(f,x)\right\}. \]  

The intervals $(x_2,x_1)$ and $(x_1,\infty)$ can be handled 
via a unified calculation, so let us complete our analysis by 
examining the intervals $(x_2,\infty)$ and $(0,x_2)$.  

\medskip 

\noindent
{\bf The Interval $\pmb{(x_2,\infty)}$:}\\
Note that $f(y+x_1)\!=\!(y+x_1)^{a_3}-A(y+x_1)^{a_2}+B$ and 
thus 
\begin{eqnarray*}
f^{(k)}(y+x_1)&=&(a_3)_k(y+x_1)^{a_3-k}-(a_2)_kA(y+x_1)^{a_2-k}. 
\end{eqnarray*} 
So, letting $\phi(y)\!:=\!f(y+x_1)$, we have  for all 
$k\!\geq\!2$: 
\scriptsize 
\begin{eqnarray} 
\label{eqn:1}
\left|\frac{\phi^{(k)}(y)}{k!\phi'(y)}\right |^{1/(k-1)} 
& = & \left |\frac{(a_3)_k(y+x_1)^{a_3-k}
-(a_2)_kA(y+x_1)^{a_2-k}}
{k!(a_3(y+x_1)^{a_3-1}-a_2A(y+x_1)^{a_2-1})} 
\right|^{1/(k-1)} \\ 
\label{eqn:2}
& = & 
\left |\frac{(a_3)_k(y+x_1)^{a_3-a_2}-(a_2)_kA}{k!a_3(y+x_1)^{k-1}
\left((y+x_1)^{a_3-a_2}-\frac{a_2A}{a_3} \right)} \right |^{1/(k-1)} \\
\label{eqn:3}
& = & \frac{1}{|y|}\left |\frac{(a_3)_k}{k!a_3}
\right|^{1/(k-1)} \left |\frac{y^{k-1}\left(
(y+x_1)^{a_3-a_2}-\frac{(a_2)_kA}{(a_3)_k}\right)} 
{(y+x_1)^{k-1}\left((y+x_1)^{a_3-a_2}-
x_1^{a_3-a_2}\right)} \right |^{1/(k-1)}, 
\end{eqnarray} 
\normalsize 
which in turn is bounded above by 
\begin{eqnarray}
\label{eqn:33}
\frac{1}{|y|}\cdot\frac{a_3-1}{2}
\left |\frac{y^{k-1}}{(y+x_1)^{k-1}}\cdot
\frac{(y+x_1)^{a_3-a_2}-\frac{(a_2)_kA}{(a_3)_k}}
{(y+x_1)^{a_3-a_2}-x_1^{a_3-a_2}}\right |^{1/(k-1)}
\end{eqnarray} 
since $\left |\frac{(a_3)_k}{k!a_3}\right |^{1/(k-1)}\le 
\frac{a_3-1}{2}$. 

Now let $z\!:=\!y/x_1$ and let 
$\rho$ be $\left.\left(\frac{(a_2)_kA}{(a_3)_k}
\right)^{1/(a_3-a_2)}\right/x_1$ or $0$ according as $k\!\leq\!a_2$ or 
not. Then for $y\!\neq\!-x_1$ (or $z\!\neq\!-1$, equivalently) we have:  
\scriptsize 
\begin{eqnarray} 
\text{\scalebox{.8}[1]{$\left |\frac{y^{k-1}}{(y+x_1)^{k-1}}\cdot
\frac{(y+x_1)^{a_3-a_2}-\frac{(a_2)_kA}{(a_3)_k}}
{(y+x_1)^{a_3-a_2}-x_1^{a_3-a_2}}\right |^{1/(k-1)}$}}
\label{eqn:4}
& = & \left |\frac{z^{k-1}}{(z+1)^{k-1}}\cdot
\frac{(z+1)^{a_3-a_2}-\rho^{a_3-a_2}}{(z+1)^{a_3-a_2}-1}
\right|^{1/(k-1)}\\
\label{eqn:5}
& = & \left |\frac{z^{k-2}}{(z+1)^{k-1}}\cdot
\frac{(z+1)^{a_3-a_2}-\rho^{a_3-a_2}}{1+(z+1)+\cdots +(z+1)^{a_3-a_2-1}}
\right|^{1/(k-1)}\\
\label{eqn:6}
 & = & \left |\frac{z}{z+1}\right|^{(k-2)/(k-1)}
\left|\frac{(z+1)^{a_3-a_2}-\rho^{a_3-a_2}}{(z+1)+(z+1)^2\cdots 
+(z+1)^{a_3-a_2}}\right|^{1/(k-1)},  
\end{eqnarray} 
\normalsize 
where the last equality follows from the identity 
$\frac{z}{(z+1)^N-1}=\frac{z+1-1}{(z+1)^N-1}=\frac{1}{1+(z+1)+
\cdots+(z+1)^{N-1}}$.  

Since $\frac{z}{z+1}\!=\!1-\frac{1}{z+1}$ is clearly an increasing function of 
$z$ for all $z\!>\!-1$ (or $y\!>\!-x_1$, equivalently), and since 
$0\!>\!-\frac{x_2}{x_1}\!>\!-1$, we must have 
for all $y\!\geq\!x_2-x_1$ (or $z\!\geq\!\frac{x_2}{x_1}-1$, equivalently) 
that 
\scriptsize 
\begin{eqnarray*}
\left |\frac{z}{z+1}\right|^{(k-2)/(k-1)}& \leq& \max\left\{
\left(\frac{x_1}{x_2}-1\right)^{(k-1)/(k-2)},1\right\}\\
& = & \max\left\{\left(\left(\frac{a_3-1}{a_2-1}\right)^{1/(a_3-a_2)}-1
\right)^{(k-1)/(k-2)},1\right\}\\ 
(\star)\hspace{2cm}\mbox{} & \leq &  
\max\left\{\frac{a_3-1}{a_2-1}-1,1\right\}.  
\end{eqnarray*} 
\normalsize
Furthermore, since $\frac{(a_2-2)_k}{(a_3-2)_k}\!\leq\!1$ for 
$2\!\leq\!k\!\leq\!a_2-2$, we must have 
$\frac{x_2}{x_1}\!\geq\!\rho$ and thus the last factor of 
Equality (\ref{eqn:6}) is bounded above by $1$. So  
by Equalities (\ref{eqn:1})--(\ref{eqn:6}) and ($\star$) we obtain that 
\[|y|\left |\frac{\phi^{(k)}(y)}{ k!\phi'(y)}\right |^{1/(k-1)}
\le \frac{1}{2}\max\{(a_3-1)(a_3-2),a_3-1\}\leq (D-1)(D-2)/2,\] 
since $D\!\geq\!3$. Thus, the contribution of the interval 
$(x_2,\infty)$ to $\bar{\alpha}(D,m)$ is no more than $(D-1)(D-2)/2$.  \qed 

\medskip 

\noindent 
{\bf The Interval $\pmb{(0,x_2)}$:} 
Similar to the last case, we have for all $x\!\in\!(0,x_2)$ that 
\begin{eqnarray*} 
\left|\frac{f^{(k)}(x)}{k!f'(x)}\right |^{1/(k-1)}
& =& \left|\frac{(a_3)_kx^{a_3-k}-(a_2)_kAx^{a_2-k}}
{k!(a_3x^{a_3-1}-a_2Ax^{a_2-1})}\right|^{1/(k-1)}\\ 
& = & \left |\frac{(a_3)_kx^{a_3-a_2}-(a_2)_kA}
{k!a_3x^{k-1}\left(x^{a_3-a_2}-\frac{a_2A}{a_3}\right)}
\right |^{1/(k-1)}\\
& = & \frac{1}{x}\left|\frac{(a_3)_k}{k!a_3}
\right|^{1/(k-1)}\left|\frac{x^{a_3-a_2}-\frac{(a_2)_kA} 
{(a_3)_k}}{x^{a_3-a_2}-x_1^{a_3-a_2}}
\right|^{1/(k-1)}\\ 
& \le & \frac{1}{ x}\cdot \frac{a_3-1}{2}
\left |\frac{x^{a_3-a_2}-\frac{(a_2)_kA}{(a_3)_k} 
}{x^{a_3-a_2}-x_1^{a_3-a_2}}\right |^{1/(k-1)}.
\end{eqnarray*} 
So, for $0\!<\!z\!<\!x_2/x_1\!=\!(\frac{a_2-1}{a_3-1})^{1/(a_3-a_2)}\!<\!1$, 
we have 
\begin{eqnarray*}
\left |\frac{x^{a_3-a_2}-\frac{(a_2A)_k}{(a_3)_k}} 
{x^{a_3-a_2}-x_1^{a_3-a_2}}\right |^{1/(k-1)}
& = & \left |\frac{z^{a_3-a_2}-\rho^{a_3-a_2}}{ z^{a_3-a_2}-1}
\right|^{1/(k-1)}\\ 
& = & \left |1-\frac{1-\rho^{a_3-a_2}}{ 1 - z^{a_3-a_2}}\right |^{1/(k-1)}\\ 
& \le & \max\left\{1,\ \left |\frac{1-\rho^{a_3-a_2}}{1 - z^{a_3-a_2}}
\right |^{1/(k-1)}\right\} \ \text{ (recalling that } \rho\!<\!1\text{)}\\ 
& \le & \max\left\{1,\ \left |\frac{1}{ 1 - z^{a_3-a_2}}
\right|^{1/(k-1)}\right\} \ \text{ (recalling that } \rho\!<\!1\text{)}  \\ 
& \le & \left |\frac{a_2-1}{a_3-a_2}\right |^{1/(k-1)}\\ 
& \le & a_2-1 \\
& \leq & a_3-2. 
\end{eqnarray*} 
Thus \[ \left |\frac{f^{(k)}(x)}{ k!f'(x)}\right |^{1/(k-1)} 
\!\le\!\frac{1}{y}\cdot \frac{(a_3-1)(a_2-1)}{2},\]
and the quantity we need to majorize on this interval is no more than 
$(D-1)(D-2)/2$. \qed 

Putting together the bounds we've found over our two preceding intervals, 
we see that our upper bound for $\bar{\alpha}(D,3)$ holds and we are done. \qed 

\section{Acknowledgements} 
The authors thank Felipe Cucker, 
John McDonald, Tom Michiels, and Steve Smale for some delightful 
conversations. We also thank two anonymous referees for useful  
suggestions. 

\bibliographystyle{amsalpha}

\end{document}